\def\R{\mathbb{R}}
\def\C{\mathbb{C}}

\def\Sp{\mathbb{S}}

\def\ni{\noindent}
\def\i{\indent}

\def\half{\frac{1}{2}}

\def\half{\frac{1}{2}}

\documentclass[10pt,draft,reqno]{amsart}
     \makeatletter
     \def\section{\@startsection{section}{1}%
     \z@{.7\linespacing\@plus\linespacing}{.5\linespacing}%
     {\bfseries
     \centering
     }}
     \def\@secnumfont{\bfseries}
     \makeatother
\theoremstyle{definition}

\theoremstyle{remark}

\numberwithin{equation}{section}
\usepackage{amsfonts}

\begin{document}

\title[Brownian Manifolds]{Brownian Manifolds, Negative Type And Geo-temporal Covariances}

\author{N. H. Bingham}
\address{N. H. Bingham: Mathematics Dept., Imperial College, London SW7 2AZ}
\email{n.bingham@ic.ac.uk}
\urladdr{http://wwwf.imperial.ac.uk/~bin06/}

\author[Aleksandar MIJATOVI\'{C}]{Aleksandar MIJATOVI\'{C}}
\address{Aleksander Mijatovi\'c: Mathematics Dept., King's College, London WC2R 2LS}
\email{aleksandar.mijatovic@kcl.ac.uk}
\urladdr{https://nms.kcl.ac.uk/mijatovic/}

\author{Tasmin L. Symons}
\address{T. L. Symons: Mathematics Dept., Imperial College, London SW7 2AZ}
\email{tls11@ic.ac.uk}



\begin{abstract}
We survey {\it Brownian manifolds} -- manifolds that can parametrise Brownian motion -- and those that cannot.  We consider covariances of space-time processes, particularly those when space is the sphere -- geo-temporal processes.  There are connections with functions of negative type. 
\end{abstract}

\maketitle

\section{Brownian manifolds and negative type}

\ni {\it Spatio-temporal processes}

\i A stochastic process, $X = \{ X_t \}$ say, is generally a mathematical model for a random phenomenon evolving with time $t$ -- a {\it temporal} process; sometimes the relevant parameter is a point in space -- a {\it spatial process}, or a {\it random field}; sometimes one needs both time and space -- a {\it spatio-temporal} process.  Our main interest here is the case when space is a sphere, which we think of as the Earth; in this case we speak of a {\it geo-temporal process}. 

\i The simplest manifold that might be used for the space variable is Euclidean space, of dimension $n$ say, ${\R}{^n}$, when one might speak of a process with multi-dimensional ($n$-dimensional) time.  The next simplest space manifold is the sphere, $\Sp$ -- the Earth, say, though we shall take the $n$-sphere ${\Sp}^n$, an $n$-dimensional manifold embedded in ${\R}^{n+1}$.  

\i The most useful class of random fields for modelling purposes is the {\it Gaussian random fields}.  In this as in any other context, the prototypical Gaussian process is {\it Brownian motion}. 

\i Processes will be real-valued, unless otherwise stated. 

\medskip

\ni {\it L\'evy's Brownian motion with multi-dimensional time} 

\i One can define Brownian motion $B = (B_t: t \in \R)$ on the line (on $(\Omega, \mathcal{F}, P)$, say) as the centred Gaussian process with {\it incremental variance} 
$$
var(B_t - B_s) = |t - s|
$$ 
and (say) $B_0 = 0$.  One can regard $B$ as a map $t \mapsto B_t$ from $\R$ to the Hilbert space $H = L^2(\Omega, \mathcal{F}, P)$, and then
$$
\Vert B_t - B_s \Vert^2 = |t - s|,
$$
the left being the incremental variance.  The covariance is then given by the inner product
$$
c(t,s) := (B_t, B_s) = \half (|t| + |s| - |t - s|),
$$
and as independence is just zero correlation in the Gaussian case, this gives independent increments as usual. 

  \i Generalising this viewpoint, Paul L\'evy \cite{Lev1} showed that one can define {\it multi-parameter Brownian motion} (or Brownian motion with multi-dimensional time) as the real-valued centred Gaussian process $B = (B_t: t \in {\R}^n)$ with incremental variance
$$
i(t,s) := E[(B_s - B_t)^2] = \Vert B_t - B_s \Vert^2 = |t - s| 
$$
(using $|.|$ for Euclidean distance); see also \cite{Lev2, Lev3} for later treatments.  \\
\i One has
$$
i(s,t) = c(s,s) +c(t,t) - 2 c(s,t),                     \eqno (i-c)
$$
and as $c(s,s) = E[B_s^2] = i(s,0)$, 
$$c(s,t) = \half (i(s,0) + i(t,0) - i(s,t)).            \eqno (c-i)
$$
Thus either of $c$, $i$ determines the other; $i$ is more convenient here. 
  
\i L\'evy also showed that Brownian motion can be defined so as to be parametrised by the sphere ${\Sp}^n$, 
in addition to ${\R}^n$ as above.  Now the incremental variance is given by the geodesic distance $d$ on the sphere (from the North Pole $o$, which plays the role of the origin above):
$$
i(s,t) = \Vert B_s - B_t \Vert^2 = d(s,t).                                       \eqno ( \ast) 
$$
Thus $\sqrt{d}(s,t) = \Vert B_t - B_s \Vert$; one calls $\sqrt{d}$ a {\it Hilbertian 
distance}. 

\i A word on terminology: our incremental variance is also known by several other names: the {\it variogram} (a term due to Matheron, arising in mining), the structure function (Yaglom), mean-squared difference (Jowett), etc.; see e.g. Cressie (\cite{Cre}, 2.3.1). 

\medskip  

\ni {\it Brownian and non-Brownian manifolds} 

\i For $M$ a manifold with geodesic distance $d$, or more generally with $(M,d)$ a metric space, one can proceed as above and call $B = (B_x: x \in M)$ a Brownian motion  {\it  parametrised} by $M$ if the $B_x$ are centred Gaussian, the incremental variance is the geodesic distance, 
$$
var(B_x - B_y) = d(x,y),
$$ 
and the finite-dimensional distributions are Gaussian (that is, linear combinations $\sum c_i B_{t_i}$ are Gaussian).  Then as before, $(\ast)$ above is satisfied with $d$ the geodesic distance on $M$.  Call such a manifold, or metric space, {\it Brownian}.  Thus Euclidean space and spheres are Brownian, by L\'evy's results above.  Further examples are given by the real or complex hyperbolic spaces, a result due to Faraut and Harzallah (\cite{Far1}, III.3, \cite{FarH}, Prop. 7.3) (and implicit in Gangolli \cite{Gan}).  By contrast, quaternionic hyperbolic spaces are not Brownian (\cite{Far1}, Cor. IV.2, or \cite{FarH}), and nor is the octonion (Cayley) projective plane. \\   
\i The question of whether a space $M$ is Brownian is thus purely geometric, as it depends on whether  a map $B$ exists satisfying $(\ast)$.  See e.g. Cartier (\cite{Car}, Th. 1 d) for this viewpoint, and for background on Gaussian Hilbert spaces, Janson \cite{Jan}. 

\medskip

\ni {\it Spaces and kernels of negative type} 

\i Call a metric space $(M,d)$ of {\it negative type} if
$$
{\sum}_{i,j = 1}^n d(x_i, x_j) u_i u_j \leq 0
$$
for all $n = 2, 3, \cdots$, all points $x_i \in M$ and all real $u_i$ with $\sum u_i = 0$ (the term {\it conditionally negative definite} is also used, reflecting the condition $\sum u_i = 0$).  Call $M$ of {\it strictly negative type} if the sum above is negative for all such $u_i$ not all zero.  Such spaces are important in a variety of contexts, and have been studied at length in the books by Blumenthal \cite{Blu} and Deza and Laurent \cite{DezL}. 

\i A kernel $k: M \times M \to {\R}_+$ is of {\it negative type} if
$$
{\sum}_{i,j = 1}^n k(x_i, x_j) u_i u_j \leq 0
$$
for all $n = 2, 3, \cdots$, all points $x_i \in M$ and all real $u_i$ with $\sum u_i = 0$, and of {\it positive type} (or positive definite) if
$$
{\sum}_{i,j = 1}^n k(x_i, x_j) u_i u_j \geq 0
$$
for all $n = 2, 3, \cdots$ and all points $x_i \in M$; similarly for {\it strictly positive type}. 

\i Covariances $c$ are of positive type.  So, incremental variances $i$ are of negative type: the first two terms on the right of $(i-c)$ contribute $0$ to the relevant summation, as $\sum u_i = 0$, so the sum is $\leq 0$ as $c$ is of positive type. 

\i For negative type on locally compact groups, see Heyer (\cite{Hey2}, Ch. 5), Berg and Forst (\cite{BerF}, II). 

\medskip

\ni {\it Schoenberg's theorems} 

\i It was shown by Schoenberg (\cite{Sch1, Sch2}) in 1937-8 that a metric space $(M,d)$ is of negative type if and only if there is a map $\phi: M \to H$ for some Hilbert space $H$ with
$$
d(x,y) = \Vert \phi(x) - \phi(y) \Vert^2.
$$
Thus, when $H = L^2(\Omega, \mathcal{F}, P)$ as before, $M$ is Brownian if and only if it is of negative type, and then Brownian motion $B$ on (parametrised by) $M$ is the map $\phi$ above.  Then $\phi: (M, \sqrt{d}) \to H$ is called the {\it Brownian embedding} (or just, embedding).  See Lyons \cite{Lyo1} for a short proof of Schoenberg's theorem. 

\i The other classical result of Schoenberg relevant here \cite{Sch2} is that a kernel $k$ is of negative type iff $e^{-t k}$ is of positive type for every $t \geq 0$.  This, of course, suggests the L\'evy-Khintchine formula, and was part of Gangolli's motivation for his theory of {\it L\'evy-Schoenberg kernels} \cite{Gan}; see also \cite{Hey2}, \cite{BerF}. 

\medskip

\ni {\it The Kazhdan property} 

\i The geometrical property of being Brownian has an algebraic interpretation in the case $M = G/K$ of a symmetric space (see Section~2 for these and other related terms).  

\i Kazhdan \cite{Kaz} defined a locally compact group to have Property $(T)$, now called the {\it Kazhdan property}, if the unit representation is isolated in the space of unitary representations.  Groups with the Kazhdan property -- {\it Kazhdan groups} -- have proved to be important in many areas; for a monograph treatment, see Bekka et al. \cite{BekHV}. (We note that a locally compact group is compact iff it is amenable and Kazhdan; see e.g. Paterson \cite{Pat} for background on amenability.) 

\i The irreducible unitary representations are in bijection with the positive
definite {\it spherical functions} (Section 2 below); the set of these is
called the {\it spherical dual}.  In the rank-one case considered below, this
can be identified with a set $\Lambda \subset \R$, 
where if $M$ is compact $\Lambda$ is a discrete set
tending to infinity, while if $M$ is Euclidean, or is real or complex
hyperbolic space, $\Lambda = [0, \infty)$.  By contrast, if $M$ is quaternionic
hyperbolic space, $\Lambda = \{ 0 \} \cup [{\lambda}_0, \infty)$, where
${\lambda}_0 > 0$ (Faraut \cite{Far1, FarH}; cf. Kostant (\cite{Kost}, 428)).  Thus $M$ is Kazhdan
in this case.  Here $$
M = G/K, \qquad G = Sp(n,1), \qquad K = Sp(n) \times Sp(1)
$$
(see e.g. (\cite{BekHV}, 3.3) for a full treatment), and $Sp(n,1)$ is Kazhdan. So too is the octonion (Cayley)
projective plane. By contrast, real and complex hyperbolic space are not Kazhdan, this being most easily seen as a consequence of Schoenberg's theorem (\cite{BekHV}, 2.11).  This explains the Faraut-Harzallah results above.

\section{Symmetric spaces, spherical functions and weights}

\ni {\it Symmetric spaces} 

\i A {\it symmetric space} (Helgason \cite{Hel1, Hel2, Hel3}, Wolf \cite{Wol2}) is a Riemannian manifold $M$ whose curvature tensor is invariant under parallel translation.  These may also be described as spaces where at each point $x$ the {\it geodesic symmetry} exists: this fixes $x$ and reverses the (direction of) geodesics through $x$, an involutive automorphism (\cite{Wol2}, Ch. 11).  Then $M$ is a Riemannian homogeneous space $M = G/K$, where $G$ is a closed subgroup of the isometry group of $M$ containing the transvections, and $K$ is the isotropy subgroup of $G$ fixing the base-point $x$.  Here $(G,K)$ is called a Riemannian symmetric pair.  The Banach algebra $L_1(K\backslash G/K)$ of (Haar) integrable functions on $G$ bi-invariant under $K$ is {\it commutative}.  Pairs with this property are called {\it Gelfand pairs}, and such Banach algebras are called {\it commutative spaces} \cite{Wol2} (a misnomer, as it is the algebra, rather than the space, which is commutative). 

\i Symmetric spaces may be split into compact, Euclidean and non-compact parts; there is a duality between the compact and non-compact cases, with the Euclidean case being self-dual (\cite{Hel1}, V.2). We confine ourselves here to symmetric spaces of {\it rank one} (\cite{Hel1}, V.6).  These are {\it two-point homogeneous spaces}; they may be classified, as spheres, Euclidean and hyperbolic spaces, of {\it constant curvature} $\kappa >0$, $\kappa =0$ and $\kappa < 0$ respectively (\cite{Hel1} Sections IX.5, X.3 p.401, and \cite{Wol1}).   

\medskip

\ni {\it Spherical functions} 

\i For harmonic analysis in this context, one needs the analogue of the classical Fourier transform in Euclidean space and the Gelfand transform for Banach algebras.  This involves spherical measures and spherical functions (\cite{Wol2}, Ch. 8) and the spherical transform (\cite{Wol2}, Ch. 9); cf. Applebaum \cite{App2}. 

For $(G,K)$ a Gelfand pair, a spherical measure $m$ is a $K$-bi-invariant multiplicative linear functional on $C_c(K\backslash G/K)$; a {\it spherical function} is a continuous function $\omega: G \to \C$ such that the measure $m_{\omega}(f):=\int_G f(x) \omega(x^{-1}) d{\mu}_G(x)$ is spherical.  The map 
$$
f \mapsto \hat f(\omega) := m_{\omega}(f)=\int_G f(x) \omega(x^{-1}) d{\mu}_G(x)
$$ 
is called the {\it spherical transform} for $(G,K)$.  The {\it positive definite} spherical functions $\phi$ on $(G,K)$ are in bijection with the irreducible unitary representations $\pi$ of $G$ with a $K$-fixed unit vector $u$ via
$$
\phi(g) = \langle u, \pi(g) u \rangle.
$$
These form the {\it spherical dual} (called $\Lambda$ in \cite{Far1}, $\Omega$ in \cite{Far3}). 

\medskip

\ni {\it Weights} 

\i When $G$ is compact, the $\pi$ here are in bijection with the {\it dominant weights}, in the sense of the Cartan-Weyl theory of weights; see e.g. Applebaum (\cite{App1}, Ch. 2), or Wolf (\cite{Wol2}, 6.3).  In the rank-one case, the dominant weights are a subset $\Lambda \subset \R$, as in Section~1; here $\Lambda$ is specified by the {\it Cartan-Helgason theorem} (\cite{Wol2}, 11.4B), (\cite{Hel3}, V.1.1, 534-538, 550), (\cite{GooW}, 549-550). 
   
\i The spherical functions satisfy an integral equation due to Harish-Chandra (see e.g. \cite{Hel1}, X), which subsumes a number of the addition formulae of classical special-function theory (see e.g. (\cite{Wat}, XI) for Bessel functions and Gegenbauer polynomials).  

\medskip

\ni {\it Compact Lie groups} 

\i Compact connected Lie groups are themselves symmetric spaces (\cite{Hel1}, IV.6).  It was shown by Baldi and Rossi \cite{BalR} that $SU(2)$ is Brownian, but that $SO(n)$ is not for $n \geq 3$ 
(the question is decided by the signs of the coefficients in the Peter-Weyl expansion).  This is despite $SU(2)$ and $SO(3)$ being locally isomorphic: $SO(3) \cong SU(2)/\{\pm e \}$; cf. \cite{Far2}, chapters 7 and 8.

\section{Geo-temporal covariances}

\ni {\it Sphere cross line} 

\i For modelling purposes in the earth sciences and climatology, one needs both a space coordinate on the sphere and a time coordinate on the line (or half-line); thus the space 
$M = \Sp \times {\R}$ (or $M = \Sp \times {\R}_+$) is needed.  The most basic process one might wish to model on $M$ is Brownian motion.  But the product can be taken in several different senses, and it turns out that the question of {\it existence} of Brownian motion depends on which kind of product we take.  Recall that by L\'evy's results of \S 1, Brownian motion exists on both $\Sp$ and $\R$ (or ${\R}_+$), since both are of negative type. 

\i First, take the product of metric spaces, under Hamming distance (``city-block metric", for those who know Manhattan), under which distances $s$ add:
$$
s := s_1 + s_2,  
$$
in the obvious notation.  From the definition of negative type, this property is preserved under such products; see e.g. \cite{Blu}, \S 3.2.  So Brownian motion on the sphere cross line exists, with the product taken in this sense.  

\i Next, one can take the product under the ordinary cartesian (or pythagorean) rule:
$$
s^2 := {s_1}^2 + {s_2}^2.
$$
Here again, Brownian motion exists.  McKean \cite{McK} gives a thorough study of the white-noise case (from which the Brownian case follows by integration), starting from the work of Chentsov \cite{Che} on white noise in this setting.  McKean's construction moves between Euclidean space ${\R}^{d+1}$ and `sphere cross half-line', $\Sp^d \times {\R}_+$. 

\i By contrast, if one takes the cartesian product of two Riemannian manifolds, distance is given by the differential cartesian rule:
$$
{ds}^2 := {ds}_1^2 + {ds}_2^2,
$$
again in the obvious notation.  It turns out that $M = \Sp \times {\R}$ is no longer of negative type -- so is no longer Brownian -- viewed as a manifold in this way.  The same holds for any product of manifolds with at least one spherical factor -- or even a factor with two pairs of antipodal points.  See \cite{HjoKM} for background and details, \cite{Lee}, \cite{Pet} for Riemannian manifolds. \\  
\i Thus Brownian motion exists on $M = \Sp \times {\R}$, regarded as a product of metric spaces in both the above senses, though not of Riemannian manifolds.  This provides a starting-point for geo-temporal modelling -- but separates the effects of space and time.  

\medskip

\ni {\it Separable and non-separable covariances} 

\i The last result, however, does not take us very far.  Real applications, e.g. to climate science, involve both spatial and temporal variation together rather than separately.  The need thus arises for a range of examples of {\it non-separable} covariances on $M = \Sp \times \R$, that can be used flexibly for modelling.  For detail here we refer to  Cressie and Huang \cite{CreH}, Gneiting (\cite{Gne1} on ${\R}^n \times \R$, \cite{Gne2} on ${\Sp}^n$, \cite{GneGG}).  For space-time modelling in general, see Kyriakidis and Journel \cite{KyrJ}, Finkenst\"adt et al. \cite{FinHI}, Porcu et al. \cite{PorM, PorMS}. For applications to global climate data see the work of Jun and Stein \cite{JunS1, JunS2}, and Jeong and Jun \cite{JeoJ}.  We note that progress was hampered in the past by lack of an adequate range of examples of covariances for modelling purposes -- to the extent that practical statisticians and climatologists felt themselves forced to use as `covariances' functions that were not even positive definite.  It was the modelling needs of climate scientists, together with the interest of the Brownian case, that prompted this study.    

\section{ Complements}

\ni {\it Testing for independence} 

\i The ideas above found a new and powerful application in statistics, in the work of Sz\'ekely, Rizzo and Bakirov \cite{SzeRB} in 2007 and Sz\'ekely and Rizzo \cite{SzeR1} in 2009.  See in particular the extensive commentary to the invited paper \cite{SzeR1}, and for further developments, \cite{SzeR2, SzeR3, SzeR4, SzeR5}; there are also applications to the theory of algorithms.  Given a bivariate sample $((X_1,Y_1), \cdots, (X_n,Y_n))$, where each coordinate has finite mean, it turns out that one can test for independence of the $X$- and $Y$-coordinates, consistently against all alternatives (again, with finite means) by test statistics involving only {\it distances} between observations.  The crux is the concept of {\it distance covariance} (equivalent to a related concept of {\it Brownian covariance}): see below. 

\newpage

\ni {\it Distance covariance} 

\i The theory of distance covariance in metric spaces has been re-worked and simplified by Russell Lyons \cite{Lyo1, Lyo2}.  It turns out that this area too belongs to distance geometry.  The crux is for the distance covariance of $(X,Y)$ to be zero iff $X$ and $Y$ are independent.  It turns out that this does not hold for general metric spaces, but does so exactly for those of {\it strong negative type}, a class that includes Euclidean spaces.  As before, embeddability into Hilbert space is crucial; other embeddings also occur (Riesz, Fourier, Crofton, Brownian).  We refer to the excellent paper \cite{Lyo1} for details.  For applications to high-dimensional data, see Kosorok \cite{Koso}. 

\i Regarding the link with Crofton's formula: see the paper by Guyan Robertson \cite{Rob}.  The Crofton formula goes back to 1885 \cite{Cro}; see Santal\'o \cite{San} for background and details. It is a precursor of the Radon transform, for which see e.g. \cite{Hel4}.

\medskip  

\ni {\it Other approaches} 

\i The first person to use white-noise integrals for L\'evy's Brownian motion was Chentsov \cite{Che}, an approach later taken up by L\'evy himself \cite{Lev3}, and McKean \cite{McK}.  For other approaches, see the work of Molchan \cite{Mol1, Mol2, Mol3, Mol4}, Noda \cite{Nod1, Nod2}, and  Takenaka, Kubo and Urakawa \cite{TakKU}. 

\medskip

\ni {\it Gaussian processes} 

\i One may construct Gaussian processes wherever one has a covariance, though these are necessarily more complicated than Brownian motion when working on a non-Brownian space.  Covariance structure is always important, but is only fully informative in the case of Gaussianity -- here, as always, a useful first approximation.  For compact symmetric spaces (such as spheres), a detailed study was given by Askey and Bingham \cite{AskB}.  It would be interesting to extend this study to the geo-temporal context. 

\medskip

\ni {\it Gaussian fields} 

\i Gaussian random fields (this term is now more common for spatial processes) have been studied by Cohen and Lifshits \cite{CohL} on spheres and hyperbolic spaces.  For Gaussian free fields -- which arise in physics (quantum field theory), but may be regarded as multi-parameter analogues of Brownian motion -- see Sheffield \cite{She}.  For an extensive survey of Gaussian random fields and their links with physics, see L\'eandre \cite{Lea}.  For contours in this context (motivated by work of Pyke on Brownian sheets), see Kendall \cite{Ken1}.

\medskip

\ni {\it Positive definite functions on spheres} 

\i This subject, at the heart of the work here, goes back to  seminal work by Schoenberg \cite{Sch3} in 1942.  It has been considered further by the first author \cite{Bin2}, Faraut \cite{Far1} and Gneiting \cite{Gne2}. 

\medskip  

\ni {\it Manifold-valued Brownian motion} 

\i Dual (in the sense of harmonic analysis -- see Gangolli \cite{Gan}) to (real-valued) Brownian motion {\it parametrised} by a manifold is Brownian motion {\it taking values in} a manifold.  Here there is a rich interplay between the geometry of the manifold and the probabilistic properties of Brownian motion on it.  See for example Grigoryan \cite{Gri}, Varopoulos et al. \cite{VarSC}, Elworthy \cite{Elw}. One of the highlights in this context is that the radial part of a Brownian motion on a Riemannian manifold is 
a semimartingale, even though the smoothness of the distance function fails along the cut-locus, see Kendall \cite{Ken2}. For a probabilistic proof of the Atiyah-Singer index theorem using Brownian motion on manifolds, see Bismut \cite{Bis}.

\medskip

\ni {\it Hypergroups} 

\i The theory of hypergroups is by now well established, but too extensive for us to discuss here.  We refer to the standard work on the subject by Bloom and Heyer \cite{BloH1}; see also \cite{BloH2, BloR1, BloR2,Hey3}.  Hypergroups make contact with the work studied here, for instance through our main example, the symmetric spaces of rank one; these have constant curvature $\kappa$.  For the spherical case $\kappa > 0$, the relevant hypergroup here is the Bingham (or Bingham-Gegenbauer) hypergroup; see \cite{Bin1}, Bloom and Heyer \cite{BloH1}, 3.4.23.  For the Euclidean case $\kappa = 0$, it is the Kingman (or Kingman-Bessel) hypergroup (\cite{BloH1}, 3.4.30).  For the hyperbolic case $\kappa < 0$, it is the hyperbolic hypergroup (\cite{BloH1}, Section 3.5.68, \cite{Zeu}). 

\medskip

\ni {\it Markov property} 

\i In one dimension, the Markov property is expressed by present time being a {\it splitting time}: past and future are {\it conditionally independent given the present}.  In higher dimensions, the geometry is more complicated.  In the plane, for example, one might have values within and without a contour conditionally independent given values on the contour.  See for example Evstigneev \cite{Evs} for background and references. 

\medskip

\ni {\it Time series} 

\i The work above is relevant to time series.  For instance, Zhou \cite{Zho} used distance correlation to study nonlinear dependence.  Prediction theory may be extended to hypergroups -- see e.g. H\"osel and Lasser \cite{HosL} -- and can be applied on spheres, using the Bingham hypergroup.  This provided the first author with a pleasing link to his recent work on prediction theory; see e.g. \cite{Bin3, Bin4, BinM}. 

\medskip

\ni {\it Fractional processes} 

\i Brownian motion is too smooth for some purposes, and may be usefully generalised to fractional Brownian motion, which has a parameter (essentially the Hurst parameter from hydrology) that governs the degree of smoothness.  Such fractional Gaussian fields have been studied in contexts related to ours by Gelbaum \cite{Gel} and Istas \cite{Ist}. 

\medskip
 
\ni {\it Higher dimensions} 

\i It is of interest to see what happens to the $n$-dimensional spheres and
hyperbolic spaces considered here as the dimension $n \to \infty$.  There has
been much progress in this area in recent decades, due largely to Olshanski,
Okounkov and Vershik.  For background and details, see several recent papers by
Jacques Faraut, e.g. \cite{Far3}; cf. the paper by Bloom and Wildberger \cite{BloW} in
this volume. 

\section{Postscript}
\i It is a pleasure to dedicate this survey to Herbert Heyer on the occasion of
his eightieth birthday.  The preface to his classic book \cite{Hey1} is preceded by
a quote from Pierre Lelong: {\it ` \ldots les probabilit\'es sur les structures
alg\'ebriques, sujet neuf et passionnant'}.  This captures well the lifelong
dedication to the subject that Herbert has always shown.  His work, his example
and his friendship have enriched the lives of us, our fellow-contributors and
many others; long may they continue to do so. 

\par\bigskip\noindent
{\bf Acknowledgment.} We are most grateful to the referee for many helpful comments and references, which have greatly improved the paper.

\bibliographystyle{amsplain}

\begin{thebibliography}{99}

\bibitem{App1} Applebaum, D.: {\sl Probability on compact Lie groups}.  Springer, 2014.
\bibitem{App2} Applebaum, D.: Convolution semigroups of probability measures on Gelfand pairs, revisited.  This volume.
\bibitem{AskB} Askey, R. A. and Bingham, N. H.: Gaussian processes on compact symmetric spaces.
{\sl Z. Wahrschein.  verw. Geb.} {\bf 37} (1976), 127 -- 143. 
\bibitem{BalR} Baldi, P. and Rossi, M.: On L\'evy's Brownian motion indexed by elements of compact groups.  {\sl Colloq. Math.} {\bf 133}.2 (2012), 227 -- 236. 
\bibitem{BekHV} Bekka, B., De la Harpe, P. and Villette,  A.:  {\sl Kazhdan's property (T)}.  Cambridge University Press, 2005. 
\bibitem{BerF} Berg, C. and Forst, G.: {\sl Potential theory on locally compact abelian groups}.  Erg. Math. {\bf 87}, Springer, 1975. 
\bibitem{Bin1} Bingham, N. H.: Random walk on spheres.  {\sl Z. Wahrschein.} {\bf 22} (1972), 169-192. 
\bibitem{Bin2}   Bingham, N. H.: Positive definite functions on spheres.  {\sl Proc. Cambridge Phil. Soc.} {\bf 73} (1973), 145--156. 
\bibitem{Bin3}  Bingham, N. H.: Szeg\H{o}'s theorem and its probabilistic descendants.  {\sl Probability Surveys} {\bf 9} (2012), 287--324. 
\bibitem{Bin4}  Bingham, N. H.:, Probability unfolding, 1965--2015.  {\sl Probability, Analysis and Number Theory} (N. H. Bingham Festschrift, ed. C. M. Goldie and A. Mijatovi\'c),
{\sl Advances in Applied Probability} {\bf 48A} (2016), 1--13. 
\bibitem{BinM} Bingham, N. H. and Missaoui, B.:  Aspects of prediction. {\sl J. Appl. Prob.} {\bf 51A} (2014), 189--201. 
\bibitem{Bis} Bismut, J. M.: Probability and geometry. {\sl Probability and analysis}, 1 -- 60,  {\sl Lecture Notes in Math.} {\bf 1206}, Springer, 1986.
\bibitem{BloH1} Bloom, W. R. and Heyer, H.: {\sl Harmonic analysis of probability measures on hypergroups}.  De Gruyter Studies in Math. {\bf 20}, Walter de Gruyter, Berlin, 1994.
\bibitem{BloH2} Bloom, W. R. and Heyer, H: Negative definite functions and convolution semigroups of probability measures on a commutative hypergroup.  {\sl Probability and Mathematical Statistics} {\bf 16} (1996), 157--176. 
\bibitem{BloR1} Bloom, W. R. and Ressel, P.: Representations of negative definite functions on polynomial hypergroups.  {\sl Arch. Math. (Basel)} {\bf 78} (2002), 318--328. 
\bibitem{BloR2} Bloom, W. R. and Ressel, P: Negative definite and Schoenberg functions on commutative hypergroups.  {\sl J. Austral. Math. Soc.} {\bf 79} (2005), 25--37. 
\bibitem{BloW} Bloom, W. R.  and Wildberger, N. J., Positive definite functions on spheres and hyperbolic spaces.  This volume.
\bibitem{Blu} Blumenthal, L. M.: {\sl Theory and applications of distance geometry}, 2nd ed., Chelsea, New York, 1970. 
\bibitem{Car} Cartier, P.: Introduction \`a l'\'etude des mouvements browniens \`a plusieurs param\`etres.  {\sl S\'em. Prob. V}, 58--75, {\sl Lecture Notes in Math.} {\bf 191}, Springer, 1971. 
\bibitem{Che} Chentsov, N. N.: L\'evy Brownian motion for several parameters and generalised white noise.  {\sl Th. Prob. Appl.} {\bf 2} (1957), 265--266. 
\bibitem{CohL} Cohen, S. and Lifshits, M. A.: Stationary Gaussian random fields on hyperbolic spaces and on Euclidean spheres.  {\sl ESAIM Probability and Statistics} {\bf 16} (2012), 165--221. 
\bibitem{Cre} Cressie, N.: {\sl Statistics for spatial data}, revised ed., Wiley, 1993.
\bibitem{CreH}  Cressie, N. and Huang, H. C.: Classes of non-separable spatio-temporal stationary covariance functions.  {\sl J. Amer. Stat. Assoc.} {\bf 94} (1999), 1330-1340. 
\bibitem{Cro} Crofton, M. W.: Probability.  {\sl Encyclopaedia Britannica}, 9th ed. {\bf 19} (1885), 768-788. 
\bibitem{DezL} Deza, M. M. and Laurent,M.: {\sl Geometry of cuts and metrics}.  Springer, 1997.
\bibitem{Elw} Elworthy, K. D.: Geometric aspects of diffusions on manifolds. {\sl \'Ecole d'\'Et\'e de Probabilit\'es de Saint-Flour XVII}, 276 -- 465, {\sl Lecture Notes in Math.} {\bf 1362}, Springer, 1988; 
\bibitem{Evs} Evstigneev, I. V.: Stochastic extremal problems and the strong Markov property of random fields. {\sl Russian Math. Surveys} {\bf 43} (1988), 1--49. 
\bibitem{Far1} Faraut, J.: Fonction brownienne sur une vari\'et\'e Riemannienne.  {\sl S\'em. Prob. VII}, 61-76, {\sl Lecture Notes in Math}.  {\bf 321}, Springer, 1973. 
\bibitem{Far2} Faraut, J.: {\sl Analysis on Lie groups}.  Cambridge Studies in Advanced Math. {\bf 110}, Cambridge Univ. Press, 2008. 
\bibitem{Far3} Faraut, J.: Asymptotics of spherical functions for large rank: An introduction.  {\sl Representation theory, complex analysis and integral geometry} (ed. B. Kr\"otz, O. Offen and E. Sayag, Biurkh\"auser/Springer, 2012), 251--275.
\bibitem{FarH} Faraut, J. \& Harzallah, K.: Distances hilbertiennes invariantes sur un espace homog\`ene.  {\sl Ann. Inst. Fourier} {\bf 24}.3 (1974), 171---283. 
\bibitem{FinHI} Finkenst\"adt, B., Held, L. and Isham, V.: {\sl Statistics of spatio-temporal systems}.  Chapman \& Hall/CRC, 2007. 
\bibitem{Gan} Gangolli, R.: Positive definite kernels on homogeneous spaces and certain stochastic processes related to L\'evy's Brownian motion of several parameters.  {\sl Ann. Inst. H. Poincar\'e B (NS)} {\bf 3} (1967), 121--226. 
\bibitem{Gel} Gelbaum, Z. A.: Fractional Brownian fields over manifolds.  {\sl Trans. Amer. Math. Soc.} {\bf 366}.9 (2014), 4781--4814. 
\bibitem{Gne1} Gneiting, T.: Non-separable, stationary covariance functions for space-time data.  {\sl J. Amer. Stat. Assoc.} {\bf 87} (2002), 590-600. 
\bibitem{Gne2} Gneiting, T.: Strictly and non-strictly positive definite functions on spheres.  {\sl Bernoulli} {\bf 19} (2013), 1327 -- 1349. 
\bibitem{GneGG} Gneiting, T., Genton, M. G. and Guttorp, P.: Geostatistical space-time models, stationarity, separability and full symmetry.  In [33], 151 -- 175. 
\bibitem{GooW} Goodman, R. and Wallach, N.: {\sl Representations and invariants of the classical groups}.  Cambridge University Press, 1998.
\bibitem{Gri} Grigor'yan, A.: Analytic and geometric background of recurrence and transience of the Brownian motion on Riemannian manifolds.  {\sl Bull. Amer. Math. Soc.} {\bf 36} (1999), 135--249. 
\bibitem{Hel1} Helgason, S.: {\sl Differential geometry and symmetric spaces}.  Academic Press, 1962. 
\bibitem{Hel2} Helgason, S.: {\sl Differential geometry, Lie groups and symmetric spaces}.  Academic Press, 1978 (2nd ed., Grad. Studies in Math. {\bf 34}, Amer. Math. Soc., 2001). 
\bibitem{Hel3} Helgason, S.: {\sl Groups and geometric analysis: Integral geometry, invariant differential operators and spherical functions}.  Academic Press, 1984 (2nd ed., Amer. Math. Soc., 2008). 
\bibitem{Hel4} Helgason, S.: {\sl Geometric analysis on symmetric spaces}.  Math. Surveys and Monographs {\bf 39}, Amer. Math. Soc., 1994. 
\bibitem{Hey1} Heyer, H.: {\sl Probability meaures on locally compact groups}.  Ergebnisse Math. {\bf 94}, Springer, 1977. 
\bibitem{Hey2} Heyer, H.: {\sl Structural aspects in the theory of probability: A primer in probabilities on algebraic-topological structures}.  World Scientific, 2004.  
\bibitem{Hey3} Heyer, H.: Positive and negative definite functions on a hypergroup and its dual.  {\sl Infinite-dimensional harmonic analysis IV}, 63--96, World Scientific, 2009. 
\bibitem{HjoKM} Hjorth, P. G., Kokkendorff. S. L. and Markvorsen, S.: Hyperbolic spaces are of strictly negative type.  {\sl Proc. Amer. Math. Soc.} {\bf 130} (2002), 175---181. 
\bibitem{HosL}  H\"osel, V. and Lasser, R.: Prediction of weakly stationary sequences on polynomial hypergroups.  {\sl Annals of Probability} {\bf 31} (2003), 93--114. 
\bibitem{Ist} Istas, J.: Multifractional Brownian fields indexed by metric spaces with distances of negative type.  {\sl ESIAM Prob. Stat.} {\bf 17} (2013), 219--223. 
\bibitem{Jan} Janson, S.: {\sl Gaussian Hilbert spaces}.  Cambridge Tracts in Math. {\bf 129}, Cambridge University Press, 1997.
\bibitem{JeoJ} Jeong, J. and Jun, M.: Covariance models on the surface of a sphere: when does it matter?  {\sl Stat.} {\bf 4} (2015), 167 -- 182.
\newpage
\bibitem{JunS1} Jun, M. and Stein, M. L.: An approach to producing space-time covariance functions on spheres. {\sl Technometrics} {\bf 49} (2007), no. 4, 468 -- 479. 
\bibitem{JunS2} Jun, M. and Stein, M. L.: Nonstationary covariance models for global data. {\sl Ann. Appl. Stat.} {\bf 2} (2008), no. 4, 1271 -- 1289.
\bibitem{Kaz} Kazhdan, D.: Connection of the dual space of a group with the structure of its closed subgroups.  {\sl Funct. Anal. Appl.} {\bf 1} (1967), 63--65. 
\bibitem{Ken1} Kendall, W. S.: Contours of Brownian processes with several-dimensional times. {\sl Z. Wahrsch.} {\bf 52} (1980), 267--276. 
\bibitem{Ken2} Kendall, W. S.: The radial part of Brownian motion on a manifold: a semimartingale property. {\sl Ann. Prob.} {\bf 15} (1987), 1491--1500. 
\bibitem{Koso} Kosorok, M. R.: On Brownian distance covariance and high-dimensional data.  {\sl Annals of Applied Statistics} {\bf 3} no. 4 (2009), 1270--1278; Correction, ibid. {\bf 7} no. 2 (2013), 1247. 
\bibitem{Kost} Kostant, B.: On the existence and irreducibility of certain series of representations.  {\sl Bull. Amer. Math. Soc.} {\bf 75} (1969), 627--642.
\bibitem{KyrJ} Kyriakidis, P. C. and Journel, A. C.: Geostatistical space-time models: a review.  {\sl Math. Geology} {\bf 31} (1999), 651 -- 684. 
\bibitem{Lea} L\'eandre, R.: The geometry of Brownian surfaces. {\sl Probability Surveys} {\bf 3} (2006), 37--88. 
\bibitem{Lee} Lee, J. M.  {\sl Riemannian manifolds: An introduction to curvature}.  Grad. Texts in Math. {\bf 176}, Springer, 1997. 
\bibitem {Lev1} L\'evy, P. {\sl Processus stochastiques et mouvement brownien}, Gauthier-Villars, Paris, 1948 (2nd ed. 1965).  
\bibitem{Lev2} L\'evy, P.: Le mouvement brownien fonction d'un point sur la sph\`ere de Riemann.  {\sl Rend. Circ. Mat. Palermo} {\bf 8} (1959), 297--310. 
\bibitem{Lev3} L\'evy, P.: Fonctions browniennes dans l'espace Euclidien et dans l'espace de Hilbert.  {\sl Research Papers in Statistics} (J. Neyman Festschrift, ed. F. N. David), Wiley, 1966, 189--223. 
\bibitem{Lyo1} Lyons, R.: Distance covariance in metric spaces.  {\sl Annals of Probability} {\bf 41} (2013), 3284--3305. 
\bibitem{Lyo2} Lyons, R.: Hyperbolic space has strong negative type.  {\sl Illinois J. Math.} {\bf 58} (2014), 1009--1013. 
\bibitem{McK} McKean, H.P.: Brownian motion with several-dimensional time. {\sl Th. Prob. Appl.} {\bf 8} (1963), 357--378. 
\bibitem{Mol1} Molchan, G. M.: On some problems concerning Brownian motion in L\'evy's sense.  {\sl Th. Prob. Appl.} {\bf 8} (1963), 682--690.
\bibitem{Mol2} Molchan, G. M.: On homogeneous random fields on symmetric spaces of rank 1.  {\sl Th. Prob. Math. Stat.} (1980), 143--168. 
\bibitem{Mol3} Molchan, G. M.: Multiparameter Brownian motion on symmetric spaces. {\sl Th. Prob. Math. Stat.} (1985), 275-286. 
\bibitem{Mol4} Molchan, G. M.: Multiparameter Brownian motion.  {\sl Th. Prob. Math. Stat.} (1988), 97--110. 
\bibitem{Nod1} Noda, A.: L\'evy's Brownian motion: total positivity structure of $M(t)$-processes and deterministic character, {\sl Nagoya Math. J.} {\bf 94} (1984), 137--164. 
\bibitem{Nod2} Noda, A.: Generalized Radon transform and L\'evy's Brownian motion, I, II.  {\sl Nagoya Math. J.} {\bf 105} (1987), 71-87, 89--107. 
\bibitem{Pat} Paterson, A. L. T.: {\sl Amenability}.  Math. Surveys and Monographs {\bf 29}, Amer. Math. Soc., 1988. 
\bibitem{Pet} Petersen, P.: {\sl Riemannian geometry}, 2nd ed.  Grad. Texts in Math. {\bf 171}, Springer, 2006. 
\bibitem{PorM} Porcu, E. and Meteu, J. (ed.): {\sl Positive definite functions: From Schoenberg to space-time challenges}.  Castellon de la Plana, 2008. 
\bibitem{PorMS} Porcu, E., Montero, J. and Schlather  M. (ed.): {\sl Advances and challenges in space-time modelling of natural events}.  Springer, 2010. 
\bibitem{Rob} Robertson, G. Crofton formula and geodesic distance in hyperbolic space.  {\sl J. Lie Theory} {\bf 8} (1998), 163--172. 
\bibitem{San} Santal\'o, L. A.: {\sl Integral geometry and geometric probability}.  Encycl. Math. Appl. {\bf 1}, Addison-Wesley, 1976. 
\bibitem{Sch1} Schoenberg, I. J.: On certain metric spaces arising from Euclidean spaces by a change of metric, and their embedding in Hilbert space.  {\sl Ann. Math.} {\bf 38} (1937), 787--793. 
\newpage
\bibitem{Sch2} Schoenberg, I. J.: Metric spaces and positive definite functions.  {\sl Trans. Amer. Math. Soc.} {\bf 44} (1938), 522-536 (reprinted in {\sl I. J. Schoenberg, Selected Papers, Vol. 1} (ed. C. de Boor), Birkh\"auser, 1988, 100--114). 
\bibitem{Sch3} Schoenberg, I. J.: Positive definite functions on spheres.  {\sl Duke Math. J. } {\bf 9} (1942), 96-108 ({\sl Selected Papers 1}, 172--185). 
\bibitem{She} Sheffield, S. Gaussian free fields for mathematicians.  {\sl Probab. Th. Rel. Fields} {\bf 139} (2007), 521-541. 
\bibitem{SzeR1} Sz\'ekely, G. J. and Rizzo, M. L.: Brownian distance covariance.  {\sl Annals of Applied Statistics} {\bf 3} no. 4 (2009), 1236--65.
\bibitem{SzeR2} Sz\'ekely G. J. and Rizzo M. L.: On the uniqueness of distance covariance.  {\sl Statistics and Probability Letters} {\bf 82} (2012), 2278--2282.
\bibitem{SzeR3} Sz\'ekely G. J. and Rizzo M. L.: The distance-correlation $t$-test of independence in higher dimensions.  {\sl J. Multivariate Anal.} {\bf 17} (2013), 193--213.  
\bibitem{SzeR4} Sz\'ekely G. J. and Rizzo M. L.: Energy statistics: A class of statistics based on distance.  {\sl J. Statist. Planning and Inference} {\bf 143} (2013), 1249--1272. 
\bibitem{SzeR5} Sz\'ekely G. J. and Rizzo M. L.: Partial distance correlation with methods for dissimilarities.  {\sl Annals of Statistics} {\bf 42} (2014), 2382--2412. 
\bibitem{SzeRB} Sz\'ekely G. J. and Rizzo M. L.: Measuring and testing dependence by correlation of distances.  {\sl Ann. Statist.} {\bf 35} (2007), 2769--2794. 
\bibitem{TakKU} Takenaka, S., Kubo, I. and Urakawa, H.: Brownian motion parameterised with metric spaces of constant curvature.  {\sl Nagoya Math. J.} {\bf 82} (1981), 131--140. 
\bibitem{VarSC} Varopoulos, N. T. , Saloff-Coste, L. and Coulhon, T.: {\sl Analysis and geometry on groups}.  Cambridge Tracts in Mathematics {\bf 100}, Cambridge University Press, 1992. 
\bibitem{Wat} Watson, G. N.: {\sl Theory of Bessel functions}, 2nd ed., Cambridge University Press, 1962 (1st ed. 1922). 
\bibitem{Wol1} Wolf, J. A.: {\sl Spaces of constant curvature}.  Amer. Math. Soc., 1967 (6th ed. 2011). 
\bibitem{Wol2} Wolf, J. A.: {\sl Harmonic analysis on commutative spaces}.  Amer. Math. Soc., 2007.
\bibitem{Zeu} Zeuner, H.: On hyperbolic hypergroups. {\sl Probability on Groups VIII}, 216 -- 224, {\sl  Lecture Notes in Math.}, {\bf 1210}, Springer, 1986.
\bibitem{Zho} Zhou Z.: Measuring nonlinear dependence in time series: a distance correlation approach.  {\sl J. Time Series Analysis} {\bf 33} (2012), 438--457. 

\end{thebibliography}

\end{document}